# REJOINDER

By S. C. Kou, Qing Zhou and Wing H. Wong

*Harvard University, Harvard University and Stanford University*

We thank the discussants for their thoughtful comments and the time they have devoted to this project. As a variety of issues have been raised, we shall present our discussion in several topics, and then address specific questions asked by particular discussants.

**1. Sampling algorithms.** The widely used state-of-the-art sampling algorithms in scientific computing include temperature-domain methods, such as parallel tempering and simulated tempering, energy-domain methods, such as multicanonical sampling and the EE sampler, and methods involving expanding the sampling/parameter space. The last group includes the Swendsen–Wang type algorithms for lattice models, as Wu and Zhu pointed out, and the group Monte Carlo method [1]. If designed properly, these sampling-space-expansion methods could be very efficient, as Wu and Zhu's example in computer vision illustrated. However, since they tend to be problem-specific, we did not compare the EE sampler with them. The comparison in the paper is mainly between the EE sampler and parallel tempering. Atchadé and Liu's comparison between the EE sampler and the multicanonical sampling thus complements our result. It has been more than 15 years since multicanonical sampling was first introduced. However, we feel that there are still some conceptual questions that remain unanswered. In particular, the key idea of multicanonical sampling is to produce a flat distribution in the energy domain. But we still do not have a simple intuitive explanation of (i) why focusing on the energy works, (ii) why a distribution flat in the energy is sought, and (iii) how such a distribution helps the sampling in the original sample space. The EE sampler, on the other hand, offers clear intuition and a visual picture: the idea is simply to "walk" on the equi-energy sets, and hence focusing on the energy directly helps avoid local trapping. In fact, the numerical results in Atchadé and Liu's comment clearly demonstrate the advantage of EE over multicanonical sampling in the 20 normal mixture example. Specifically, their Table 1 shows that in terms of estimating the probabilities of visiting each mode, the EE sampler is about two to three times more efficient. We think that estimating the probability of visiting individual modes provides a more sensitive measure of the performance, the reason being that even if a sampler misses two or three modes in each run, the sample average of the first and second moments could still be quite good; for example, missing one mode in the far lower left can be offset by missing one mode in the far upper right in the sample average of the first moment, and missing one faraway mode can be offset by disproportionately visiting much more frequently another faraway mode in the sample average of the second moment, and so on. Nevertheless, we agree with Atchadé and Liu that more studies (e.g., on the benchmark phase transition problems in the Ising and Potts models) are needed to reach a firmer conclusion.

**2. Implementing the EE sampler for scientific computations.** The EE sampler is a flexible and all-purpose algorithm for scientific computing. For a given problem, it could be adapted in several ways.

First, we suggested in the paper that as a good initial start the energy and temperature ladders could be both assigned through a geometric progression. It is conceivable that for a complicated problem alternative assignments might work better, as Minary and Levitt's off-lattice protein folding example illustrated. A good assignment makes the acceptance rates of the EE jump comparably across the different chains, say all greater than 70%. This can be achieved by a small pilot run of the algorithm, which can be incorporated into an automatic self-tuning implementation.

Second, the energy ladder and temperature ladder can be decoupled in the sense that they do not need to always obey $(H_{i+1} - H_i)/T_i \approx c$. For example, for discrete problems such as the lattice phase transition models and the lattice protein folding models, one could take each discrete energy level itself as an energy ring, while keeping the temperatures as a monotone increasing sequence. In this case an EE jump is always accepted, since it always moves between states with the same energy level.

Third, the EE sampler can be implemented in a serial fashion as Wu and Zhu commented. One could start the algorithm from $X^{(K)}$, run for a predetermined number of iterations, completely stop it and move on to $X^{(K-1)}$, run it, completely stop, move on to $X^{(K-2)}$, and so on. This serial implementation offers the advantage of saving computer memory in that one only needs to record the states visited in the chain immediately preceding the current one. The downside is that it will not provide the users the option to online monitor and control (e.g., determine to stop) the algorithm; instead, one has to prespecify a fixed number of iterations to run. In the illustrative multimodal distribution in the paper and the example we include in this rejoinder in Section 4, we indeed utilized the serial implementation since the number of iterations for each chain was prespecified.

Fourth, the EE sampler constructs energy rings to record the footsteps of high-order chains. The fact that a computer's memory is always finite might appear to limit the number of iterations that the EE sampler can be run. But as Minary and Levitt pointed out, this seeming limitation can be readily solved by first putting an upper bound (subject to computer memory) on the energy ring size; once this upper bound is reached a new sample can be allocated to a specific energy ring by replacing a randomly chosen element in the ring. Minary and Levitt's example involving a rough one-dimensional energy landscape provides a clear demonstration.

Fifth, the key ingredient of the EE sampler is the equi-energy move, a global move that compensates for the local exploration. It is worth emphasizing that the local moves can adopt not only the Metropolis–Hastings type moves, but also Gibbs moves, hybrid Monte Carlo moves as in Minary and Levitt's example, and even moves applied in molecular dynamic simulations, as long as the moves provide good explorations of the local structure.

Sixth, the equi-energy move jumps from one state to another within the same energy ring. As Wu and Zhu commented, it is possible to conduct moves across different energy rings. It has pros and cons, however. It might allow the global jump a larger range, and at the same time it might also lead to a low move acceptance rate, especially if the energy of the current state differs much from that of the proposal jump state. The latter difficulty is controlled in the equi-energy jump of the EE sampler, since it always moves within an energy ring, where the states all have similar energy levels. One way to enhance the global jump range and rein in the move acceptance rate is to put a probability on each energy ring in the jump step. Suppose the current state is in ring $D_j$. One can put a distribution on the ring index so that the current ring $D_j$ has the highest probability to be chosen, and the neighboring rings $D_{j-1}$ and $D_{j+1}$ have probabilities less than that of $D_j$ to be chosen, and rings $D_{j-2}$ and $D_{j+2}$ have even smaller probabilities to be chosen, and so on. Once a ring is chosen, the target state is proposed uniformly from it.

**3. Theoretical issues.** We thank Atchadé and Liu for providing a more probabilistic derivation of the convergence of the EE sampler that complements the one we gave in the paper. While these results assure the long-run correctness of the sampler, we agree, however, with Wu and Zhu that investigating the convergence speed is theoretically more challenging and interesting, as it is the rate of convergence that separates different sampling algorithms. So far the empirical evidence supports the EE sampler's promise, but definitive theoretical results must await future studies.

In addition to facilitating the empirically observed fast convergence, another advantage offered by the idea of working on the equi-energy sets is that it allows efficient estimation by utilizing all the samples from *all* the chains on an energy-by-energy basis (as discussed in Section 5 of the paper). We thus believe that the alternative estimation strategy proposed by Chen and Kim is very inefficient, because it essentially wastes all the samples in the chains other than the target one. To make the comparison transparent, suppose we want to estimate the probability of a rare event under the target distribution $P_{\pi_0}(X \in A)$. Chen and Kim's formula would give

$$\hat{P} = \frac{1}{n} \sum_{i=1}^{n} \sum_{j=0}^{K} w_j \mathbb{1}(X_i^{(0)} \in A \cap D_j).$$

But since $P_{\pi_0}(X \in A)$ is small, say less than $10^{-10}$, there is essentially no sample falling into $A$ in the chain $X^{(0)}$, and correspondingly $\hat{P}$ would be way off no matter how cleverly $w_j$ is constructed. The fact that the high-order chains $X^{(j)}$ could well have samples in the set $A$ (due to the flatness of $\pi_j$) does not help at all in Chen and Kim's strategy. But in the EE estimation method such high-order-chain samples are all employed. The tail probability estimation presented in Section 5 and Table 4 illustrates the point. The reason that the EE estimation method is much more efficient in this scenario is due to the well-known fact that in order to accurately estimate a rare event probability importance sampling has to be used and the fact that the EE strategy automatically incorporates importance sampling in its construction. We also want to point out that rare event estimation is an important problem in science and engineering; examples include calculating surface tension in phase transition in physics, evaluating earthquake probability in geology, assessing the chance of bankruptcy in insurance or bond payment default in finance, estimating the potentiality of traffic jams in telecommunication, and so on.

**4. Replies to individual discussants.** We now focus on some of the individual points raised. Minary and Levitt's discussion has been covered in Sections 1 and 2 of this rejoinder, as was Wu and Zhu's in Sections 1 to 3; we are sorry that space does not permit us to discuss their contributions further.

Atchadé and Liu questioned the derivation of (5) of the paper. This equation, we think, arises directly from the induction assumption, and does not use any assumption on $X^{(i+1)}$ explicitly or implicitly. We appreciate their more probabilistic proof of the convergence theorem.

Chen and Kim asked about the length of the burn-in period in the examples. In these examples the burn-in period consists of 10% to 30% of the samples. We note that this period should be problem-dependent. A rugged high-dimensional energy landscape requires longer burn-in than a smooth low-dimensional one. There is no one-size-fits-all formula.

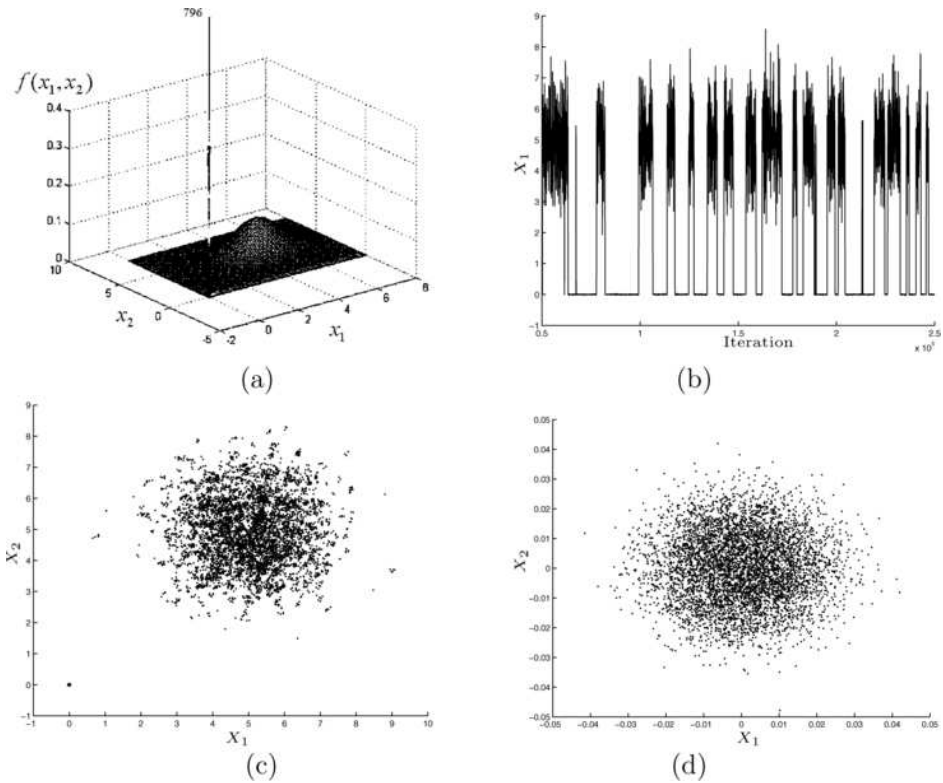

Fig. 1. *The artificial needle-in-the-haystack example.* (a) *The density function of the target distribution.* (b) *The sample path of* $X_1^{(0)}$ *from a typical run of the EE sampler.* (c) *The samples generated at both modes. Note the mode at the origin.* (d) *The samples generated near the mode at the origin.*

In the discussion Chen and Kim appeared to suggest that the Gibbs sampler is preferred in high-dimensional problems. But our experience with the Gibbs sampler tells a different story. Though simple to implement, in many cases the Gibbs sampler can be trapped by a local mode or by a strong correlation between the coordinates—the very problems that the modern state-of-the-art algorithms are trying to tackle.

We next consider the needle-in-the-haystack example raised in Chen and Kim's discussion, in which the variances of the normal mixture distribution differ dramatically. Figure 1(a) shows the density function of this example. We implemented the EE sampler using four chains (i.e., $K = 3$) and 200,000 iterations per chain after a burn-in period of 50,000 iterations. Following the energy ladder setting used in Chen and Kim, we set $H_1 = 3.13$; the other energy levels were set between $H_1$ and $H_{\min} + 100$ ($= 93$) in a geometric progression: $H_1 = 3.13, H_2 = 8.3, H_3 = 26.8$. The MH proposals were specified as $N_2(\mathbf{X}_n^{(i)}, \tau_i^2 T_i I_2)$, where $T_i$ ($i = 0, \ldots, K$) is the temperature of the $i$th chain. We set $\tau_i = 1$ for $i > 0$ and $\tau_0 = 0.05$. The probability of equi-energy jump $p_{\text{ee}} = 0.3$. With all the above parameters fixed in our simulation, we tested the EE sampler with different highest temperatures $T_K$, whereas the remaining temperatures were evenly distributed on the log-scale between $T_K$ and $T_0 = 1$. We tried $T_K =$10, 20, 30, 50 and 100; with each parameter setting the EE sampler was performed independently 100 times. From the target chain $\mathbf{X}^{(0)}$ we calculated

$$\hat{P} = \frac{1}{n} \sum_{i=1}^{n} \mathbb{1}(\sqrt{(X_{i1}^{(0)})^2 + (X_{i2}^{(0)})^2} < 0.05),$$

the probability of visiting the mode at the origin. From the summary statistics in Table 1, we see that (i) the performance of EE is quite stable with an MSE between 0.04 and 0.06 for different temperature ladders; (ii) more than 98% of the times EE did jump between the two modes. In order to assess the performance of EE on this problem, we also applied PT under exactly the same settings including the numbers of chains and iterations, the temperature ladders and the exchange probability ($p_{\text{ex}} = p_{\text{ee}} = 0.3$). It turns out that with all the different temperature ladders PT never outperformed even the worst performance of EE ($T_K = 10$) in MSE (Table 1). From the best performance of the two methods, that is, EE with $T_K = 30$ and PT with $T_K = 20$, one sees that (i) the MSE of EE is about 54% of that of PT; (ii) the spread of the estimated probability is smaller for EE than for PT [see the standard deviation and $(5\%, 95\%)$ quantiles]. We selected a typical run of EE in the sense that the frequency of jump between the two modes of this run is approximately the same as the average frequency, and we plotted the samples in Figure 1. The chain mixed well in each mode and the cross-mode jump is acceptable. Even in this artificially created extreme example of a needle in the haystack the performance of EE is still quite satisfactory with only four chains ($K = 3$). It is worth emphasizing that we did not even fine-tune the energy or temperature ladders—they are simply set by a geometric progression.

But we do want to point out that one can always cook up extreme examples to defeat any sampling algorithm. For instance, one can hide two needles miles apart in a high-dimensional space, and no sampling algorithm is immune to this type of extreme example. In fact in Chen and Kim's example, if we ran EE with only 50,000 iterations (after the burn-in period) with $T_K = 30$, the resulting MSE increased to 0.136 and 36% of the times EE missed the needle completely (Table 1).

TABLE 1
*Summary statistics of EE and PT for the needle-in-the-haystack example*

|  | $E(\hat{P})$ | std$(\hat{P})$ | 5% | 95% | MSE | # Jump | # Miss |
|---|---|---|---|---|---|---|---|
| EE$(N=200, T_K=10)$ | 0.3740 | 0.2119 | 0.0289 | 0.7020 | 0.0603 | 36.61 | 1 |
| EE$(N=200, T_K=20)$ | 0.4298 | 0.2048 | 0.0556 | 0.7492 | 0.0464 | 40.35 | 2 |
| EE$(N=200, T_K=30)$ | 0.4567 | 0.1973 | 0.1188 | 0.7440 | 0.0404 | 43.14 | 0 |
| EE$(N=200, T_K=50)$ | 0.3958 | 0.2172 | 0.0223 | 0.6939 | 0.0576 | 39.16 | 2 |
| EE$(N=200, T_K=100)$ | 0.4396 | 0.2122 | 0.0986 | 0.7762 | 0.0482 | 39.19 | 0 |
| EE$(N=50, T_K=30)$ | 0.3077 | 0.3163 | 0 | 0.8149 | 0.1361 | 6.83 | 36 |
| PT$(N=200, T_K=10)$ | 0.4241 | 0.2971 | 0 | 0.9276 | 0.0932 | 364.07 | 7 |
| PT$(N=200, T_K=20)$ | 0.4437 | 0.2692 | 0.0000 | 0.9476 | 0.0749 | 157.18 | 4 |
| PT$(N=200, T_K=30)$ | 0.4664 | 0.3181 | 0 | 0.9979 | 0.1013 | 104.20 | 6 |
| PT$(N=200, T_K=50)$ | 0.4793 | 0.3093 | 0 | 0.9204 | 0.0951 | 63.47 | 6 |
| PT$(N=200, T_K=100)$ | 0.4291 | 0.2972 | 0 | 0.9772 | 0.0925 | 36.02 | 7 |

Tabulated are the mean, standard deviation, 5% and 95% quantiles, and MSE of $\hat{P}$ in 100 independent runs. Also reported here are the average number of jumps between the two modes and the total number of runs in which the sampler missed the mode at the origin. $N$ is the number of iterations for each chain in units of 1000 after the burn-in period.

**5. Concluding remarks.** We thank all the discussants for their insightful contributions. We appreciate the efforts of the Editor and the Associate Editor for putting up such a platform for exchanging ideas. We hope that the readers will enjoy as much as we did reading these comments and thinking about various scientific, statistical and computational issues raised.

S. C. KOU
Q. ZHOU
DEPARTMENT OF STATISTICS
HARVARD UNIVERSITY
SCIENCE CENTER
CAMBRIDGE, MASSACHUSETTS 02138
USA
E-MAIL: kou@stat.harvard.edu
    zhou@stat.harvard.edu

W. H. WONG
DEPARTMENT OF STATISTICS
SEQUOIA HALL
STANFORD UNIVERSITY
STANFORD, CALIFORNIA 94305-4065
USA
E-MAIL: whwong@stanford.edu